\documentclass[letterpaper]{article}
\usepackage{amssymb}
\usepackage{amsmath}
\usepackage{amsfonts}
\usepackage{amscd}
\usepackage{latexsym}
\usepackage{graphicx}
\usepackage{amsthm,enumerate}
\usepackage{cases}
\usepackage[pdftex,bookmarks,colorlinks,breaklinks=true]{hyperref}
\usepackage{color}
\input xy
\xyoption{all}

\newcounter{example}
\renewcommand{\theexample}{\arabic{example}}

\numberwithin{equation}{section}

\newbox\mybox
\def\arrover#1{\mathrel{
       \setbox\mybox=\hbox spread 1.4em
              {\hfil$\scriptstyle#1\vphantom{g}$\hfil}
       \vbox{\offinterlineskip\copy\mybox
             \hbox to\wd\mybox{\rightarrowfill}}}}

\def\ontoover#1{\mathrel{
       \setbox\mybox=\hbox spread 1.4em{\hfil$\scriptstyle#1$\hfil}
       \vbox{\offinterlineskip\copy\mybox
             \hbox to\wd\mybox{\rightarrowfill\hskip-2.8mm
                               $\rightarrow$}}}}

\begin{document}

\title{EPH-classifications \\ in Geometry, Algebra, Analysis and Arithmetic}%
\author{Arash Rastegar}%


\maketitle
\begin{abstract}
Trichotomy of Elliptic-Parabolic-Hyperbolic appears in many
different areas of mathematics. All of these are named after
the very first example of trichotomy, which is formed by
ellipses, parabolas, and hyperbolas as conic sections.
We try to understand if these classifications are justified and if
similar mathematical phenomena is shared among different cases
EPH-classification is used.
\end{abstract}
\section*{Introduction}

At first glance EPH-classification is nothing deep. A discriminant of a quadratic form being
negative, zero or positive should not be that important for people to make so much fuss about it.
Geometric similarity between the zeros of quadratic forms in affine space, could be a motivation for naming them differently.
But why is it that the elliptic, parabolic, hyperbolic names are used so widely in mathematics, in places where
it seems irrelevant to the shape of hyperbola, parabola and ellipse? The question is, does there exist a paradigm behind these
namings.

Let us consider ellipses, parabolas and hyperbolas as conic sections. Parabolas are measure zero, degenerate against hyperbolas and ellipses.
But the theory of parabolas is reacher than the hyperbolic and elliptic theories. Although hyperbolas and ellipses tend to
parabolas and many properties of parabolas are limits of properties of hyperbolas or ellipses, but there are features of parabolas which
do not have counter parts in the elliptic or hyperbolic worlds. For example, every two parabolas are similar. This is not true for ellipses or hyperbolas.

Some of these features are shared by the EPH-trichotomy of Euclidean and non-Euclidean geometries. Hyperbolic geometry and elliptic geometry both tend to Euclidean geometry which is called parabolic by a simple scale of curvature. The concept of line in Euclidean space is the limit of hyperbolic lines and elliptic lines. The limiting behavior of Euclidean geometry does not limit the richness of the structure of Euclidean geometry.
Descartian coordinates is a feature of Euclidean geometry which is not a limit of similar structures in hyperbolic or elliptic worlds.
Also, the concept of similarity only makes sense in the elliptic case. Two triangles with equal angles are equal in hyperbolic geometry and elliptic geometry. Despite these similarities with conic section one can not point to any appearance of ellipses, parabolas, or hyperbolas in
elliptic, parabolic and hyperbolic geometries. This is exactly the point which makes the study of EPH-classification interesting.

Here are different sections in this paper, where EPH-classifications are discussed:
\\
1. EPH-classifications in geometry

1.1. Euclidean and non-Euclidean geometries

1.2. Manifolds of constant sectional curvature

1.3. Hyperbolic, Parabolic and Elliptic dynamical systems

1.4. Totally geodesic 2-dimensional foliations on 4-manifolds

1.5. Normal affine surfaces with $\mathbb{C}^*$-actions
\\
2. EPH-classifications in algebra

2.1. Elements of $SL_2(\mathbb{R})$

2.2. Action of $SL_2(\mathbb{R})$ on complex numbers, double numbers, and dual numbers
\\
3. EPH-classifications in analysis

3.1. Riemann's uniformization theorem

3.2. Partial differential equations

3.3. Petersson's elliptic, parabolic, hyperbolic expansions of modular forms
\\
4. EPH-classifications in arithmetic

4.1. Rational points on algebraic curves

4.2. Hyperbolic algebraic varieties

5. Categorization of EPH phenomena
\\
6. Examples of EPH abuse
\\
7. EPH controversies
\section{EPH-classifications in geometry}

"Euclidean and non-Euclidean geometries" are named elliptic, parabolic and hyperbolic by Klein.
"Manifolds of constant sectional curvature" are named elliptic, parabolic and hyperbolic after
"Euclidean and non-Euclidean geometries".
"Totally geodesic 2-dimensional foliations on 4-manifolds" are related to "Manifolds of constant sectional curvature".

Although "Hyperbolic, Parabolic and Elliptic dynamical systems" are very different in nature, they still fit into the geometric paradigm
because of the geometric nature of dynamical systems.
"Normal affine surfaces with $\mathbb{C}^*$-actions" are named elliptic, parabolic and hyperbolic after
"Hyperbolic, Parabolic and Elliptic dynamical systems".
Therefore, the five geometric subsections seem to be of two different origins
and not coming from the original paradigm of elliptic, parabolic, hyperbolic as conic sections.

\subsection{Euclidean and non-Euclidean geometries}

The name non-Euclidean was used by Gauss to describe
a system of geometry which differs from Euclid's in
its properties of parallelism (Coxeter 1998).
Spherical geometry was not historically considered to be
non-Euclidean in nature, as it can be embedded in a 3-dimensional Euclidean space.
Following Ptolemi, muslims pioneered in understanding the spherical geometry.
The subject was unified in 1871 by Klein, who gave the names
parabolic, hyperbolic, and elliptic to the respective systems
of Euclid, Bolyai-Lobacevskii, and Riemann-Shlafli (Klein 1921).
Felix Klein is usually given credit for being the first
to give a complete model of a non-Euclidean geometry.
Roger Penrose notes that it was Eugenio Beltrami
who first discovered both the projective and conformal
models of the hyperbolic plane (Penrose 2005). Klein built his model by suitably adapting
Arthur Cayley�s metric for the projective plane. Klein went on,
in a systematic way, to describe nine types of two-dimensional geometries. Yaglom calls
these geometries Cayley-Klein geometries (McRae 2007).

Yaglom gave conformal models for these geometries, extending what had been done for both the projective
and hyperbolic planes. Each type of geometry is homogeneous and can be determined by two
real constants $\kappa_1$ and $\kappa_2$ (Yaglom 1979). The metric structure is called Elliptic, Parabolic, Hyperbolic, respectively according to if $\kappa_1 > 0$, $\kappa_1= 0$ or $\kappa_1< 0$. Conformal structure is called Elliptic, Parabolic, Hyperbolic, respectively according to if
$\kappa_2 > 0$, $\kappa_2= 0$ or $\kappa_2< 0$. 

\subsection{Manifolds of constant sectional curvature}

A hyperbolic $n$-manifold is a complete Riemannian $n$-manifold of constant sectional curvature $-1$.
Every complete, connected, simply-connected manifold of constant negative curvature $-1$ is isometric to the real hyperbolic space $\mathbb{H}^n$. As a result, the universal cover of any closed manifold $M$ of constant negative curvature $-1$ is $\mathbb{H}^n$. Thus, every such $M$ can be written as $\mathbb{H}^n/ \Gamma$ where $\Gamma$ is a torsion-free discrete group of isometries on $\mathbb{H}^n$. That is, $\Gamma$ is a lattice in $SO^{+}(1,n)$.
Its thick-thin decomposition has a thin part consisting of tubular neighborhoods of closed geodesics and ends which are the product of a Euclidean $(n-1)$-manifold and the closed half-ray. The manifold is of finite volume if and only if its thick part is compact.
For $n>2$ the hyperbolic structure on a finite volume hyperbolic n-manifold is unique by Mostow rigidity and so geometric invariants are in fact topological invariants.

The geometries of constant curvature can be classified into the following three cases:
Elliptic geometry (constant positive sectional curvature),
Euclidean geometry (constant vanishing sectional curvature),
Hyperbolic geometry (constant negative sectional curvature).
The universal cover of a manifold of constant sectional curvature is one of the model spaces:
sphere, Euclidean space, and hyperbolic upper half-space.

\subsection{Hyperbolic, Parabolic and Elliptic dynamical systems}

In general terms, a smooth dynamical system is called hyperbolic if the tangent space
over the asymptotic part of the phase space splits into two complementary directions, one
which is contracted and the other which is expanded under the action of the system. In the
classical, so-called uniformly hyperbolic case, the asymptotic part of the phase space is
embodied by the limit set and, most crucially, one requires the expansion and contraction
rates to be uniform. Uniformly hyperbolic systems are now fairly well understood. They
may exhibit very complex behavior which, nevertheless, admits a very precise description.
Moreover, uniform hyperbolicity is the main ingredient for characterizing structural
stability of a dynamical system. Over the years the notion of hyperbolicity was broadened
(non-uniform hyperbolicity) and relaxed (partial hyperbolicity, dominated splitting)
to encompass a much larger class of systems, and has become a paradigm for complex
dynamcial evolution (Arauju-Viana 2008)

The theory of uniformly hyperbolic dynamical systems was initiated in the 1960�s
(though its roots stretch far back into the 19th century) by S. Smale, his students and collaborators,
in the west, and D. Anosov, Ya. Sinai, V. Arnold, in the former Soviet Union.
It came to encompass a detailed description of a large class of systems, often with very
complex evolution. Moreover, it provided a very precise characterization of structurally
stable dynamics, which was one of its original main goals.

Let $f : M \to M$ be a diffeomorphism on a manifold $M$. A compact invariant set $\Lambda\subseteq M$ is a
hyperbolic set for $f$ if the tangent bundle over $L$ admits a decomposition
$$
T_{\Lambda}M = E^u\oplus E^s,
$$
invariant under the derivative and such that $\parallel Df^{-1} | E^u\parallel < \lambda$ and $\parallel Df | E^s\parallel < \lambda$ for some
constant $\lambda < 1$ and some choice of a Riemannian metric on the manifold. When it exists,
such a decomposition is necessarily unique and continuous. We call $E^s$ the stable bundle
and $E^u$ the unstable bundle of $f$ on the set $\Lambda$.
The definition of hyperbolicity for an invariant set of a smooth flow containing no equilibria
is similar, except that one asks for an invariant decomposition $T_{\Lambda}M = E^u\oplus E^0\oplus E^s$,
where $E^u$ and $E^s$ are as before and $E^0$ is a line bundle tangent to the flow lines. An invariant
set that contains equilibria is hyperbolic if and only it consists of a finite number of points,
all of them hyperbolic equilibria.

When the eigenvalues are on the unit circle and complex, the dynamics is called elliptic.
The elliptic paradigm revolves around two features at the opposite
end of the orbit complexity scale from the exponential behavior captured by the hyperbolic
paradigm. The first and most important is a remarkable persistence for fairly general
classes of conservative dynamical systems of stable behavior (in certain parts of the phase
space), which can be modeled on a translation of a torus. The other is somewhat less
precise. It can be roughly described as the appearance of exceptionally precise simultaneous
return of many orbits close to their initial positions. In this case no identifiable
complete set of models is available but certain typical features of both topological and
measure-theoretical behavior can be identified. The interaction between the properties of
the linearized and nonlinear systems is more subtle than in the hyperbolic case.
Both conceptually and technically, elliptic dynamics is related to hard analysis to
a greater extent than hyperbolic dynamics, where geometric and probabilistic ideas and
methods are very prominent. This is one of the reasons for the comparatively small role
elliptic dynamics plays in the field of dynamical systems. (Hasselblot and Katok 2002)

The two preceding concepts presented
an overview of the two extremes among the widespread types of phenomena in differentiable
dynamics. These may be termed �stable� (elliptic) and �random� (hyperbolic and
partially hyperbolic) motions.
What remains is a grey middle ground of subexponential behavior with zero entropy
that cannot be covered by the elliptic paradigms of stability and fast periodic approximation.
At present, there is no way to attempt a classification of the characteristic phenomena
for this kind of behavior. This is to a large extent due to the problem that the linearization
is in general not well structured and not sufficiently representative of the nonlinear
behavior. Furthermore, the outer boundary of the usefulness of the elliptic paradigm is
not sufficiently well-defined. However, there is a type of behavior
within this intermediate class that is modeled well enough by the infinitesimal and local
�shear� orbit structure of unipotent linear maps. It is this kind of behavior that we call
parabolic.(Hasselblot and Katok 2002)
\subsection{Totally geodesic 2-dimensional foliations on 4-manifolds}

In 1986 Cairns and Ghys describe the behavior of 2-dimensional totally geodesic foliations on
compact Riemannian 4-manofolds. A foliation $\mathfrak{F}$ is "geodesible" if there exists a Riemannian metric that makes
$\mathfrak{F}$ totally geodesic. The problem of chracterizing geodesible foliations has been essentially solved in the one-dimensional case by Sullivan. The codimension one case is so rigid that one can give a complete classification.The basic point is that in codimension 1, the distribution
$\mathfrak{F}^{\perp}$, orthogonal to $\mathfrak{F}$, is evidently integrable. In arbitrary codimension complicacies arise and the simplest case these can be well treated in case of 2-dimensional foliations on 4-manifolds.(Cairns and Ghys 1986)

Assuming that the foliation $\mathfrak{F}$ and manifold $M$ are oriented and $C^{\infty}$. In general 2-dimensional foliation of arbitrary codimension are treated. Cairns and Ghys prove the following results: Let $\mathfrak{F}$ be a 2 dimensional geodesible foliation on a compact manifold $M$, then there exists a Riemannian metric $g$ on $M$ for which $\mathfrak{F}$ is totally geodesicand such that the curvature of the leaves is the same constant $K$ equal to $+1$, $0$ and $-1$.These foliations are called elliptic ($K=1$), parabolic ($K=0$) and hyperbolic ($K=-1$) which are treated in the following statements: $\mathfrak{F}$ is elliptic geodesible if and only if the leaves of $\mathfrak{F}$.
are the fibers of a fibration of $M$ by spheres $S^2$. For a parabolic geodesible 2-dimensional foliation on a compact 4-manifold $M$ we know that,
there exists an Abelian covering $\widehat{M}$ of $M$ such that the lift $\widehat{\mathfrak{F}}$ of $\mathfrak{F}$
to $\widehat{M}$ can be defined by a locally free action of $\mathbb{R}^2$. But if the leaves are of constant negative curvature on a 4-manifold, the orthogonal distribution $\mathfrak{F}^{\perp}$ is necessarily integrable and there are two possibilities:

(1) either $\mathfrak{F}$ is defined by a suspension of a representation of the fundamental group of some surface of genus greater than one in the
group of diffeomorphisms of the sphere $S^2$, or

(2) the universal covering space $\widetilde{M}$ of $M$ is diffeomorphic to $\mathbb{R}^4$, in such a way that the leaves  of $\mathfrak{F}$ are covered by $\mathbb{R}^2\times \{ *\}$ and those of  $\mathfrak{F}^{\perp}$ by $\{ *\}\times\mathbb{R}^2$.

Using these reuslts Cairns and Ghys deduce that, if there exists a geodesic foliation $\mathfrak{F}$ on a compact simply connected 4-manifold $M$,
then $M$ is one of the two fibration by spheres $S^2$ over $S^2$ and the leaves of $\mathfrak{F}$ are the fibers of this foliation. If $\mathfrak{F}$
is a 2-dimensional geodesic folition on a compact 4-manifold with negative Euler characteristic, then $M$ is a fibration by spheres $S^2$ over a compact surface. Moreover one can choose this fibration in such a way that its fibers are either everywhere tangent or everywhere transverse to $\mathfrak{F}$.

\subsection{Normal affine surfaces with $\mathbb{C}^*$-actions}

A $\mathbb{C}^*$-action on a normal affine surface is called elliptic if it has a unique fixed point
which belongs to the closure of every 1-dimensional orbit, parabolic if the set of its
fixed points is 1-dimensional, and hyperbolic if has only a finite number of fixed
points, and these fixed points are of hyperbolic type, that is each one of them belongs
to the closure of exactly two 1-dimensional orbits.

In the elliptic case, the complement $V^*$ of the unique fixed point in $V$ is fibered
by the 1-dimensional orbits over a projective curve $C$. In the other two cases is
fibered over an affine curve $C$, and this fibration is invariant under the $\mathbb{C}^*$-action.

Vice versa, given a smooth curve $C$ and a $\mathbb{Q}$-divisor $D$ on $C$, the Dolgachev-
Pinkham-Demazure construction provides a normal affine surface $V=V_{C,D}$ with a
$\mathbb{C}^*$-action such that is just the algebraic quotient of $V^*$ or of $V$ , respectively. This
surface $V$ is of elliptic type if $C$ is projective and of parabolic type if $C$ is affine.
(Flenner and Zeidenberg 2003).

\section{EPH-classifications is algebra}

We shall justify below why "Elements of $SL_2(\mathbb{R})$" are named elliptic, parabolic, and hyperbolic after the original paradigm of conic sections. "Actions of $SL_2(\mathbb{R})$ on complex numbers, double numbers, and dual numbers" are also related to the geometry of conic sections. Therefore, the subsections of EPH-classifications in algebra all belong to a single paradigm.

\subsection{Elements of $SL_2(\mathbb{R})$}

A key part of the study of M\"{o}bius transformations is their classification into elliptic, parabolic and hyperbolic transformations. The classification can be summarized in terms of the fixed points of the transformation, as follows:

Parabolic: The transformation has precisely one fixed point in the Riemann sphere.

Hyperbolic: There are exactly two fixed points, one of which is attractive, and one of which is repelling.

Elliptic: There are exactly two fixed points, both of which are neutral.

Now we turn to studying the connection between M\"{o}bius transformations and conic sections. The transformations we consider are those that preserve the unit disk . It is straightforward to see that all such M\"{o}bius transformations correspond to matrices in $SU(1,1)$.
Given such a matrix , the M\"{o}bius transform  maps the unit disk into itself in such a way that it maps arcs orthogonal to the unit circle to arcs orthogonal to the unit circle. These transforms are in fact the isometries of the Poincar\'{e} disk. We can interpret them as mappings of the hyperboloid  by using natural mapping  defined from the hyperboloid model to Poincar\'{e} disk.
This way, we end up with a homomorphism $h:SU(1,1)\to SO(2,1)$. This mapping is not injective. In fact, it is a 2-to-1 mapping.

We have now shown how M\"{o}bius transformations are mapped to  matrices in $SO(2,1)$. The classification of M\"{o}bius transformations is determined by the properties of these  matrices. These are matrices $A$ such that $A^T JA=J$ where $J=diag(1,1,-1)$.
First, we must show that each such matrix a has an eigenvector  with eigenvalue 1. It has three eigenvalues, counting multiplicity. Suppose $\lambda$  is an eigenvalue (possibly complex), and let $x$ be an eigenvector for matrix $A\in SO(2,1)$. Then $Ax=\lambda x$. Since $A^T JA=J$, we have
$Jx=A^T JAx=\lambda A^T Jx$. Then $Jx$ is an eigenvector of $A^T$ with eigenvalue $1/\lambda$. Since $A$ and $A^T$ have the same eigenvalues, we conclude that $1/\lambda$ is also an eigenvalue of $A$. Thus, if $\lambda\neq 1$ , then the third eigenvalue must be 1, in order for the determinant of $A$ (which is always equal to the product of all the eigenvalues) to be 1. This proves $A$ has an eigenvalue equal to 1.

Also, if $A$ has more than one independent eigenvector with eigenvalue 1, then it must be the identity matrix. Then for non-identity matrices , let $v$ be a unit eigenvector with eigenvalue 1, i.e. such that $Av=v$. All the planes $\{v^T Jx=c\}$ are invariant under $A$ , since
$v^T JAx=v^T(A^T)^{-1}Jx=(A^{-1}v)^TJx=v^T Jx$.

This is a family of parallel planes all invariant under $A$.
Then we say that $A$ is elliptic, parabolic or hyperbolic if and only if the corresponding cross-sections of the standard right circular cone are
ellipses, parabolas, or hyperbolas, respectively.

\subsection{Action of $SL_2(\mathbb{R})$ on complex numbers, double numbers, and dual numbers}

The $SL_2(\mathbb{R})$ action by M\"{o}bius transformations is usually considered
as a map of complex numbers $z = x+iy$,
$i^2 = -1$. Moreover, this action defines a map from
the upper half-plane to itself.
However there is no need to be restricted to
the traditional route of complex numbers only.
Less-known double and dual numbers
also have the form $z = x + iy$ but different assumptions
on the imaginary unit $i$: $i^2 = 0$ or $i^2 = 1$
correspondingly. Although the arithmetic of dual
and double numbers is different from the complex
ones, e.g., they have divisors of zero, we are still
able to define their transforms by M\"{o}bius transformations in most cases.
Three possible values -1, 0, and 1 of $\sigma :=i^2$
will be referred to here as elliptic, parabolic,
and hyperbolic cases respectively. (V. Kisil 2007)

To understand the M\"{o}bius action in all EPH cases, we
use the Iwasawa decomposition of $SL_2(\mathbb{R}) = ANK$ into three one-dimensional subgroups
$A, N, K$:
$$
\left(
  \begin{array}{cc}
    a & b \\
    c & d \\
  \end{array}
\right)
=
\left(
  \begin{array}{cc}
    k & 0 \\
    0 & k^{-1} \\
  \end{array}
\right)
\left(
  \begin{array}{cc}
    1 & w \\
    0 & 1 \\
  \end{array}
\right)
\left(
  \begin{array}{cc}
    Cos \theta & -Sin\theta \\
    Sin\theta & Cos\theta \\
  \end{array}
\right).
$$
Subgroups $A$ and $N$ act in (1) irrespective of the
value of $\sigma$: $A$ makes a dilation by $\alpha^2$, i.e., $z \rightarrow \alpha^2 z$,
and $N$ shifts points to left by $\nu$, i.e. $z \rightarrow z + \nu$.

By contrast, the action of the third matrix
from the subgroup $K$ sharply depends on $\sigma$. In the elliptic, parabolic and hyperbolic
cases $K$-orbits are circles, parabolas and
(equilateral) hyperbolas correspondingly.

Now, we explore the isometric action of
$G = SL_2(\mathbb{R})$ on the homogeneous spaces
$G/H$ where $H$ is one dimensional subgroup of $G$. There are only three such subgroups
up to conjugacy:
$$
K = \{
\left(
  \begin{array}{cc}
    \cos \theta & -\sin \theta \\
    \sin \theta & \cos \theta \\
  \end{array}
\right)
|0\leq\theta <2\pi\}
$$
$$
N'= \{
\left(
  \begin{array}{cc}
    1 & 0 \\
    t & 1 \\
  \end{array}
\right)|t\in \mathbb{R}
\}
$$
$$
A'=\{
\left(
  \begin{array}{cc}
    \cosh \alpha & \sinh \alpha \\
    \sinh \alpha & \cosh \alpha \\
  \end{array}
\right)|\alpha\in\mathbb{R}
\}
$$
These three subgroups give rise to elliptic, parabolic and hyperbolic geometries
(abbreviated EPH). The name comes about from the shape of the equidistant
orbits which are ellipses, parabolas, hyperbolas respectively. Thus the Lobachevsky
geometry is elliptic (not hyperbolic) in this terminology.

EPH are 2-dimensional Riemannian, non-Riemannian and pseudo-Riemannian
geometries on the upper half-plane (or on the unit disc). Non-Riemannian geometries,
are a growing field with geometries like Finsler gaining
more influence. The Minkowski geometry is formalised in a sector of a flat plane
by means of double numbers.
Subgroups $H$ fix the imaginary unit $i$ under the above action and thus
are known as EPH rotations (around $i$). Consider a distance function invariant
under the $SL_2(\mathbb{R})$ action. Then the orbits of $H$ will be equidistant points from
$i$, giving some indication on what the distance function should be. But this does
not determine the distance entirely since there is freedom in assigning values to the
orbits. Review a well-known standard definition of distance $d : X \times X \to \mathbb{R}^+$ with:

(1) $d(x, y) \geq 0$,

(2) $d(x, y) = 0$ iff $x = y$,

(3) $d(x, y) = d(y, x)$,

(4) $d(x, y) \leq d(x, z) + d(z, y)$,
for all $x, y, z \in X$.

Although adequate in many cases, the defined concept does
not cover all interesting distance functions. For example, in the
hugely important Minkowski space-time the reverse of the triangle inequality holds.

Recall the established procedure of constructing geodesics in Riemannian geometry
(two-dimensional case):

(1) Define the metric of the space: $Edu^2 + Fdudv + Gdv^2$.

(2) Define length for a curve $\Gamma$ as:
$$
length(\Gamma) =\int_{\Gamma}(Edu^2 + Fdudv + Gdv^2)^{\frac{1}{2}}
$$

(3) Then geodesics will be defined as the curves which give a stationary point
for length.

(4) Lastly the distance between two points is the length of a geodesic joining
those two points.

In this respect, to obtain the $SL_2(\mathbb{R})$ invariant distance we require the invariant metric.
The invariant metric in EPH cases is
$$
ds^2 =\frac{du^2 -\sigma dv^2}{v^2}
$$
where $\sigma = -1, 0, 1$ respectively.
For an arbitrary curve $\Gamma$:
$$
length(\Gamma) = \int_{\Gamma}\frac{(du^2 -\sigma dv^2)^{\frac{1}{2}}}{v^2}
$$
In the two non-degenerate cases (elliptic and hyperbolic) to find the geodesics
is straightforward, it is now the case of solving the Euler-Lagrange equations and
hence finding the minimum or the maximum respectively. The Euler-Lagrange
equations for the above metric take the form:
$
\frac{d}{dt}\left(\frac{\dot{\gamma_1}}{y^2}\right)=0
$
and
$
\frac{d}{dt}\left(\frac{\dot{\gamma_2}}{y^2}\right)=\frac{\dot{\gamma_1}-\sigma\dot{\gamma_2}}{y^3}
$
where $\gamma$ is a smooth curve
$\gamma(T ) = (\gamma_1(T ),\gamma_2(T ))$ and $T \in (a, b)$.
For $\sigma = -1$ the solution is well-known: semicircles orthogonal to the real axes or
vertical lines. Equations of the ones passing though $i$ are:
$$
(x^2 + y^2) \sin 2t - 2x \cos 2t - \sin 2t = 0
$$
where $t \in \mathbb{R}$. And the distance function is then:
$$
d(z,w) = \sinh^{-1}\frac{|z - w|}{2\sqrt{\Im z\Im w}}
$$
where $\Im z$ is the imaginary part of $z$.

In the hyperbolic case when $\sigma = 1$ there are two families of solutions, one space-like,
one time-like:
$$
x^2 - y^2 - 2tx + 1 = 0
$$
$$
(x^2 - y^2) \sinh 2t - 2x \cosh 2t + sinh 2t = 0
$$
with $t \in \mathbb{R}$. Those are again to the real axes. And the distance
functions are:
$$
d(z,w) =
2 \sin^{-1} \frac{\sqrt{\Re (z-w)^2-\Im (z-w)^2}}{2\sqrt{\Re z \Im w}}
$$
when space-like;
$$
d(z,w) =
2 \sinh
^{-1} \frac{\sqrt{\Re (z-w)^2-\Im (z-w)^2}}{2\sqrt{\Re z \Im w}}
$$
when time-like, where $\Re z$ and $\Im (z - w)$ are the real and imaginary part of $z$.

In the parabolic framework the only solutions are vertical line, again orthogonal to the real axes.
They indeed minimise the distance between two points $w_1, w_2$ since the geodesic is up the line
$x = \Re w_1$ through infinity and down $x = \Re w_2$. Any points on the same vertical
lines have distance zero, so $d(w_1,w_2) = 0$ for all $w_1, w_2$ which is a very degenerate
function.(A. Kisil 2009)
\section{EPH-classifications in analysis}

"Riemann's uniformization theorem" belong to the paradigm of "Manifold of constant curvature".
"Partial differential equations" could be called elliptic, parabolic and hyperbolic according to the eigenvalues of a matrix being
negative or positive which seems to be a new paradigm different from the above.
"Petersson's elliptic, parabolic, hyperbolic expansions of modular forms" come from the same paradigm of "Elements of $SL_2(\mathbb{R})"$ which
seems to be the original paradigm of conic sections.

\subsection{Riemann's uniformization theorem}

 The uniformization theorem says that every simply connected Riemann surface is conformally equivalent to one of the three domains: the open unit disk, the complex plane, or the Riemann sphere. In particular it admits a Riemannian metric of constant curvature. This classifies Riemannian surfaces as elliptic (positively curved � rather, admitting a constant positively curved metric), parabolic (flat), and hyperbolic (negatively curved) according to their universal cover.
 The uniformization theorem is a generalization of the Riemann mapping theorem from proper simply connected open subsets of the plane to arbitrary simply connected Riemann surfaces.

Felix Klein (1883) and Henri Poincar\`{e} (1882) conjectured the uniformization theorem for (the Riemann surfaces of) algebraic curves. Henri Poincar\`{e} (1883) extended this to arbitrary multivalued analytic functions and gave informal aguments in its favor. The first rigorous proofs of the general uniformization theorem were given by Poincar\'{e} (1907) and Paul Koebe (1907). Paul Koebe later gave several more proofs and generalizations.

On an oriented surface, a Riemannian metric naturally induces an almost complex structure as follows: For a tangent vector $v$ we define $J(v)$ as the vector of the same length which is orthogonal to $v$ and such that $(v, J(v))$ is positively oriented. On surfaces any almost complex structure is integrable, so this turns the given surface into a Riemann surface.

From this, a classification of metrizable surfaces follows. A connected metrizable surface is a quotient of one of the following by a free action of a discrete subgroup of an isometry group:

1.the sphere (curvature $+1$)

2.the Euclidean plane (curvature $0$)

3.the hyperbolic plane (curvature $-1$).

The first case includes all surfaces with positive Euler characteristic: the sphere and the real projective plane. The second includes all surfaces with vanishing Euler characteristic: the Euclidean plane, cylinder, M\"{o}bius strip, torus, and Klein bottle. The third case covers all surfaces with negative Euler characteristic: almost all surfaces are hyperbolic. For closed surfaces, this classification is consistent with the Gauss-Bonnet theorem, which implies that for a closed surface with constant curvature, the sign of that curvature must match the sign of the Euler characteristic.

The positive/flat/negative classification corresponds in algebraic geometry to Kodaira dimension $-1,0,1$ of the corresponding complex algebraic curve.

\subsection{Partial differential equations}

The general second-order PDE in two independent variables has the form analogous to the equation for a conic section.
More precisely, replacing $\partial /\partial x$ by $x$, and likewise for other variables (formally this is done by a Fourier transform), converts a constant-coefficient PDE into a polynomial of the same degree, with the top degree (a homogeneous polynomial, here a quadratic form) being most significant for the classification. Just as one classifies conic sections and quadratic forms into parabolic, hyperbolic, and elliptic based on the discriminant , the same can be done for a second-order PDE at a given point. However, the discriminant in a PDE is formulated slightly differently, due to the change in convention of the $xy$ term.

1. In case of negative discriminant, solutions of elliptic PDEs are as smooth as the coefficients allow, within the interior of the region where the equation and solutions are defined.

2. In case of zero discriminant, equations that are parabolic at every point can be transformed into a form analogous to the heat equation by a change of independent variables. Solutions smooth out as the transformed time variable increases.

3. In case of positive discriminant, hyperbolic equations retain any discontinuities of functions or derivatives in the initial data.
If there are n independent variables, the classification of a general linear partial differential equation of second order depends upon the signature of the eigenvalues of the coefficient matrix.

1.Elliptic: The eigenvalues are all positive or all negative.

2.Parabolic : The eigenvalues are all positive or all negative, save one that is zero.

3.Hyperbolic: There is only one negative eigenvalue and all the rest are positive, or there is only one positive eigenvalue and all the rest are negative.

4.Ultra-hyperbolic: There is more than one positive eigenvalue and more than one negative eigenvalue, and there are no zero eigenvalues. There is only limited theory for ultra-hyperbolic equations (Courant and Hilbert, 1962).

For higher order PDEs, there are some nice descriptions of hyperbolicity for higher order equations and systems. There is generally no description for parabolicity. And elliptic, parabolic, and hyperbolic do not come close to exhausting the possible PDEs that can be written down beyond second order scalar equations.

To scratch at the surface of this problem, we need to journey back to 1926, when the definitions of elliptic, parabolic, and hyperbolic PDEs are given by Jacques Hadamard. In his study of scalar linear partial differential equations of second order (the work has since been compiled and published as Lectures on Cauchy's problem in linear partial differential equations by Dover publications in 1953), Hadamard made the following definitions. (As an aside, it is also in those lectures that Hadamard made the first modern definition of well-posedness of the Cauchy problem of a PDE.)

Given a linear partial differential operator of second order with real coefficients, its principal part can be represented by a (symmetric) matrix of coefficients $a_{ij}\partial_{ij}$  . The operator is then said to be

�Elliptic if $a_{ij}$   is positive definite or negative definite
�Parabolic if $a_{ij}$   is positive or negative semi-definite, and admits precisely one 0 eigenvalue.
�Hyperbolic if $a_{ij}$   is indefinite, but is non-degenerate; that is, for every vector $v_j$  , there exists some vector $w_i$   such that
$w_i a_{ij} v_j \neq 0$ . In other words, it has no zero eigenvalue.
Already you see that the classification is incomplete: any operator with principle part having nullity $>1$ , or having nullity $=1$  but with indefinite sign is not classified. Furthermore, the definition of hyperbolic is different from the modern one. Indeed, Hadamard made the additional definition that the operator is

�Normal hyperbolic if $a_{ij}$   is hyperbolic, and furthermore all but one of its eigenvalues have the same sign.
The wave operator is a normal hyperbolic operator. Whereas what now-a-days we call the ultrahyperbolic operators (in order to distinguish them from hyperbolic ones) were nominally considered to be hyperbolic by Hadamard's standards. Hadamard was able to show that linear, normal, hyperbolic PDEs admit well-posed initial value (Cauchy) problems for data prescribed on hypersurfaces over which tangent space the restriction of $a_{ij}$   is elliptic.

\subsection{Petersson's elliptic, parabolic, hyperbolic expansions of modular forms}

Let $\Gamma\in PSL_2(\mathbb{R})$ be a Fuchsian group of the first kind acting on the upper half plane $\mathbb{H}$.
We write $x + iy = z \in \mathbb{H}$ and set $d\mu(z)$ to be the $SL_2(\mathbb{R})$-invariant hyperbolic volume form
$dxdy/y^2$. Assume the volume of the quotient space $ \mathbb{H} /\Gamma$ is equal to $V < \infty�$. Let $S_k(\Gamma)$ be the
space of holomorphic weight k cusp forms for $\Gamma$. This is the vector space of holomorphic
functions $f$ on $\mathbb{H}$ which decay rapidly in each cusp of $\Gamma$ and satisfy the transformation
property
$$
\frac{f(\gamma z)}{j(\gamma ,z)^k}-f(z)=0
\forall\gamma\in \Gamma
$$
with $j(\left(
          \begin{array}{cc}
            a & b \\
            c & d \\
          \end{array}
        \right),z)=cz+d$ for $\gamma=\left(
          \begin{array}{cc}
            a & b \\
            c & d \\
          \end{array}
        \right)$.
        We do not assume that $\Gamma$ has cusps. If there are
none then we may ignore the rapid decay condition for $S_k$. We assume throughout that
the multiplier system is trivial and that, unless otherwise stated, $4 \leq k �\in 2\mathbb{Z}$. Using the
notation $(f|_k\gamma)(z) := f(\gamma z)/j(\gamma, z)^k$, extended to all 
$\mathbb{C}[PSL_2(\mathbb{R})]$ or $\mathbb{C}[SL_2(\mathbb{R})]$ 
by linearity,the above can be written more simply as $f|_k(\gamma -1) = 0$.

 The identity in $\Gamma$ is $I = �\left( \begin{array}{cc}
            1 & 0 \\
            0 & 1 \\
          \end{array}\right).$ The remaining elements may be partitioned into three
sets: the parabolic, hyperbolic and elliptic elements. These correspond to translations, dilations
and rotations, respectively, in $\mathbb{H}$. As is well-known, the action for parabolic
elements leads to a Fourier expansion of f associated to each cusp of $\Gamma$. The parabolic
Fourier coefficients that arise often contain a great deal of number-theoretic information.
Much less well-known are Petersson�s hyperbolic and elliptic Fourier expansions. (Imamoghlu and O'Sullivan 2008)

Parabolic expansions: We say that $\gamma$
 is a parabolic element of $\Gamma$ if its trace, $tr(\gamma)$, has
absolute value 2. Then $\gamma$
 fixes one point in $R �\cup {�\infty}$. These parabolic fixed points are
the cusps of $\Gamma$. If $\mathfrak{a}$ is a cusp of $\Gamma$ then the subgroup, $\Gamma_{\mathfrak{a}}$, of all elements in $\Gamma$ that fix
$\mathfrak{a}$ is isomorphic to $\mathbb{Z}$. Thus $\Gamma_{\mathfrak{a}} = <\gamma_{\mathfrak{a}}>$ for a parabolic generator
$\gamma_{\mathfrak{a}}\in \Gamma$. There exists a scaling
matrix $\sigma_{\mathfrak{a}} �\in SL_2(\mathbb{R})$ so that $\sigma_{\mathfrak{a}}�\infty = \mathfrak{a}$ and
$$
\sigma_{\mathfrak{a}}\Gamma_{\mathfrak{a}}\sigma_{\mathfrak{a}}=\{ \pm \left(
                                                                         \begin{array}{cc}
                                                                           1 & m \\
                                                                           0 & 1 \\
                                                                         \end{array}
                                                                       \right)\mid m\in \mathbb{Z}\}
$$
The matrix $\sigma_{\mathfrak{a}}$ is unique up to multiplication on the right by any $\left(
                                                                                      \begin{array}{cc}
                                                                                        1 & x \\
                                                                                        0 & 1 \\
                                                                                      \end{array}
                                                                                    \right)
$ with $x \in \mathbb{R}$. We label
this group as $\Gamma_{\infty}$. A natural fundamental domain for $\Gamma_{\infty}\setminus \mathbb{H}$ is the set $\mathbb{F}_{\infty}$ of all
$z \in \mathbb{H}$
with $0 \leq Re(z) < 1$. The image of this set under $\sigma_{\mathfrak{a}}$ will be a fundamental domain for $\mathbb{H}/\Gamma_{\mathfrak{a}}$.

We next define an operator $A$ that converts functions with a particular parabolic invariance
into functions with invariance as $z \to z + 1$. Similar, though slightly more elaborate,
operators will do the same for functions with hyperbolic and elliptic invariance as follows.
For any function $f$ with $f|_k\gamma_{\mathfrak{a}}=f$, define $A_{\mathfrak{a}}f :=(f|_k \sigma_{\mathfrak{a}})$. Then
$(A_{\mathfrak{a}}f)(z+1)=(A_{\mathfrak{a}}f)(z)$.
It follows that $f$ in $S_k$ implies $(A_{\mathfrak{a}}f) (z)$ has period 1 and is holomorphic on $\mathbb{H}$. It
consequently has a Fourier expansion
$$
(A_{\mathfrak{a}}f) (z) =\sum_{m\in \mathbb{Z}} b_{\mathfrak{a}}(m)\exp ^{2\pi imz}.
$$
The rapid decay condition at the cusp $\mathfrak{a}$ in the definitions of $S_k$
is then
$$
(A_{\mathfrak{a}}f) (z) = (f|_k\sigma_{\mathfrak{a}})(z)\ll_{\mathfrak{a}} \exp ^{-cy}
$$
as $y \to \infty$ uniformly in $x$ for some constant $c > 0$. This must hold at each of the cusps $\mathfrak{a}$. It
is equivalent to $f|_k\sigma_{\mathfrak{a}}$ only having terms with $m > 1$ in the above expansion. Now, we define
for $f$ in $S_k$ , the parabolic expansion of $f$ at $\mathfrak{a}$ is
$$
(f|_k\sigma_{\mathfrak{a}}) (z) =\sum_{m\in \mathbb{N}} b_{\mathfrak{a}}(m)\exp ^{2\pi imz}.
$$

Hyperbolic expansions: An element $\gamma$
 of $\Gamma$ is hyperbolic if $|tr(\gamma)| > 2$. Denote the set of all such elements $Hyp(\Gamma)$. Let $\eta = (\eta_1, \eta_2)$ be a
hyperbolic pair of points in $\mathbb{R} �\cup \{�\infty\}$ for $\Gamma$. By this we mean that there exists one element
of $Hyp(\Gamma)$ that fixes each of $\eta_1, \eta_2$. The set of all such $\gamma$
 is a group which we label $\Gamma_{\eta}$. As in the parabolic case this group is isomorphic to $\mathbb{Z}$,
$\Gamma_{\eta}=<\gamma_{\eta}>$. There exists a scaling matrix $\sigma_{\eta} �\in SL_2(\mathfrak{R})$ such that
$\sigma_{\eta}0 = \eta_1, \sigma_{\eta}�\infty = \eta_2$ and
$$
\sigma_{\eta}^{-1}\gamma_{\eta}\sigma_{\eta}=\pm \left(
                                                   \begin{array}{cc}
                                                     0 & \xi \\
                                                     \xi & 0 \\
                                                   \end{array}
                                                 \right)
$$
for $\xi \in \mathbb{R}$. This scaling matrix $\sigma_{\eta}$ is unique up to multiplication on the right by any
$\left(
   \begin{array}{cc}
     x & 0 \\
     0 & x^{-1} \\
   \end{array}
 \right)
$
with $x \in \mathbb{R}$. Replacing the generator
$\gamma_{\eta}$ by $\gamma_{\eta}^{-1}$
if necessary we may assume $\xi^2 > 1$. Let
$$
\mathbb{F}_{\eta} := {z \in  \mathbb{H} : 1 \leq |z| < \xi^2}.
$$
Then it is easy to see that $\sigma_{\eta}\mathbb{F}_{\eta}$  is a fundamental domain for $\gamma_{\eta}\mid \mathbb{H}$.
Now, for any function $f$ with $f|_k\gamma_{\eta}= f$, let
$$
(A_{\eta}f) (z) := \xi^{kz}(f|_k\sigma_{\eta})(\xi^{2z}).
$$
Then $(A_{\eta}f) (z + 1) = (A_{\eta}f) (z)$.

If $f \in S_k$ then $A_{\eta}f$ has period 1 and hence a Fourier expansion:
$$
(A_{\eta}f) (z) =\sum_{m\in \mathbb{Z}} b_{\eta}(m)\exp ^{2\pi imz}.
$$
Put $w = \xi^{2z}$ so that $\exp ^{2\pi iz = w^{\pi i}/ log \xi}$. Then $f \in S_k$ must have the following
expansion: The hyperbolic expansion of $f \in S_k$ at $\eta$ is
$$
(f|_k\sigma_{\eta}) (w) =\sum_{m\in \mathbb{Z}} b_{\eta}(m)w ^{-k/2+\pi im/log \xi}.
$$

Elliptic expansions: If $z_0 = x + iy$ in $\mathbb{H}$ is fixed by a non-identity element of $\Gamma$ it is called an elliptic point of
$\Gamma$. Such group elements necessarily have traces with absolute value less than 2 (and are
called elliptic elements). Let $\Gamma_{z_0} \subset \Gamma$ be the subgroup of all elements fixing $z_0$. This is a cyclic group of finite order
$N > 1$. Let $\varepsilon\in\Gamma $ be a generator of $\Gamma_{z_0}$ . There exists $\sigma_{z_0} \in GL(2,\mathbb{C})$ so that $\sigma_{z_0} 0 = z_0, \sigma_{z_0} \infty =\overline{z_0}$. To be explicit,
we take
$$
\sigma_{z_0} =\frac{1}{2i\beta}
\left(
\begin{array}{cc}
-\overline{z_0} & z_0 \\
-1 & 1 \\
\end{array}
\right)
,\sigma_{z_0}^{-1} =\left(
\begin{array}{cc}
 1 & -z_0 \\
  1 & -\overline{z_0} \\
  \end{array}
  \right).
 $$
Note that $\sigma_{z_0}^{-1}$ maps the upper half plane $\mathbb{H}$ homeomorphically to the open unit disc
$\mathbb{D} \subset \mathbb{C}$ centered at the origin. For any $w \in \mathbb{H}$ a calculation shows
$(\sigma_{z_0}^{-1}\varepsilon\sigma_{z_0})w=\zeta^2 w$ with
$\zeta = j(\varepsilon, z_0)$ and
$$
\sigma_{z_0}^{-1}\varepsilon\sigma_{z_0}=
\left(
  \begin{array}{cc}
    \zeta & 0 \\
    0 & \zeta \\
  \end{array}
\right).
$$
Hence $\zeta$ is a primitive 2N-th root of unity: $\zeta = \exp ^{2\pi im/(2N)}$ for some $m$ with $(m, 2N) = 1$.
There exists $m'\in  \mathbb{N}$ so that $m'm\equiv 1 mod 2N$ and $\zeta^{m'} = \exp ^{\pi i/N}$. So, replacing $\varepsilon$ by $\varepsilon^{m'}$ if
necessary, we may assume $\zeta = \exp ^{\pi i/N}$. Let $\mathbb{F}_{z_0}$ equal the central sector covering $1/N$-th of the
disc and chosen with angle $\theta$ satisfying $-\pi \mid N \leq \theta -\pi \leq \pi \mid N$.
Then $\sigma_{z_0} (\mathbb{F}_{z_0})$ is a convenient fundamental domain for $\Gamma_{z_0}\mid \mathbb{H}$. Also note that there exists
$C(z_0, \Gamma) > 0$ such that
$|z| < C(z_0, \Gamma)$ for all $z \in \sigma_{z_0} (\mathbb{F}_{z_0})$.
In other words the fundamental domain we have chosen is contained in a bounded region
of $\mathbb{H}$.

Since any $f \in S_k(\Gamma)$ is holomorphic at $z = z_0$ we see that $f(\sigma_{z_0}w)$ is holomorphic at
$w = 0$ and has a Taylor series $f(z)=\sum_n a_{z_0}(n)w^n$.
Therefore we get the simple expansion
$$
f(z) =\sum_{n\in \mathbb{N}_0} a_{z_0}(n)(\sigma_{z_0}^{-1}z)^n.
$$
More useful for our purposes is the slightly different elliptic expansion due to Petersson.
For $f, g : \mathbb{H} \to \mathbb{C}$ define
$$
(A_{z_0}f) (z) := \zeta^{kz}(f|_k\sigma_{z_0})(\zeta^{2z}),
$$
$$
A_{z_0}^{-1} g := (B_{z_0}g) |_k(\sigma_{z_0}^{-1})
$$
for $(B_{z_0}g) (z) := z^{-k/2}g
\left(\frac{N \log(z)}{2\pi i}
\right)$
We have
$A_{z_0}A_{z_0}^{-1} f = A_{z_0}^{-1}A_{z_0} f=f$ and
and
$$
(f|_k\varepsilon)(z) = f(z) \Rightarrow (A_{z_0}f) (z + 1) = (A_{z_0}f) (z),
$$
$$
g(z + 1) = g(z) \Rightarrow
(A_{z_0}^{-1} g)|_k\varepsilon=A_{z_0}^{-1} g
$$
Note that the matrices $\sigma_{z_0}$ and $\sigma_{z_0}^{-1}$
have determinants $1/(2i\beta)$ and $2i\beta$ respectively. In
this case it is convenient to normalize the stroke operator $|$ and define
$$
(f|_k \gamma) (z) :=\frac{det(\gamma)^{k/2}f(\gamma z)}{j(\gamma ,z)^k}
$$
Obviously this agrees with our previous definition when
$\gamma\in SL2(R)$.

Let $f\in S_k$ then $A_{z_0}f$ has period 1 and a Fourier expansion
$$
(A_{z_0}f) (z) =\sum_{m\in \mathbb{Z}} b_{z_0}(m)\exp ^{2\pi imz}.
$$
Put $w = \zeta^{2z} = \exp ^2\pi iz/N$ so that $\exp^{2\pi iz} = w^N$ and
$$
(f|_k \sigma_{z_0}) (w) =\sum_{m\in \mathbb{Z}} b_{z_0}(m) w^{Nm-k/2}.
$$
Since $(f|_k \sigma_{z_0})$ is holomorphic at $w = 0$ we must have non-negative powers of $w$ in the
above expansion. Thus any $f\in S_k$ satisfies the following:
The elliptic expansion of $f$ in $S_k$ at $z_0$ is
$$
(f|_k \sigma_{z_0}) (z)=\sum_{m\in \mathbb{Z}}^{Nm-k/2\geq 0} b_{z_0}(m) z^{Nm-k/2}.
$$

\section{EPH-classifications in arithmetic}

"Rational points on algebraic curves" match with "Manifolds of constant curvature" but do not seems to be of geometric nature.
We can believe them to be of the same paradigm only if the same arithmetic phenomena is shared by the geometric paradigm.
"Hyperbolic algebraic varieties" in arithmetic have no elliptic or parabolic counterparts.

\subsection{Rational points on algebraic curves}

A polynomial relation $f(x, y) = 0$ in two variables defines a curve $C_0$. If the coefficients
of the polynomial are rational numbers then one can ask for solutions of the equation
$f(x, y) = 0$ with $x, y \in \mathbb{Q}$, in other words for rational points on the curve. The set of
all such points is denoted $C_0(\mathbb{Q})$. If we consider a non-singular projective model $C$ of the
curve then topologically $C$ is classified by its genus, and we call this the genus of $C_0$ also.
Note that $C_0(\mathbb{Q})$ and $C(\mathbb{Q})$ are either both finite or both infinite. Mordell conjectured,
and in 1983 Faltings proved, that if the genus of $C_0$ is greater than or equal to two, then $C_0(\mathbb{Q})$ is finite.

The case of genus zero curves is much easier and
was treated in detail by Hilbert and Hurwitz. They explicitly reduce to the cases
of linear and quadratic equations. The former case is easy and the latter is resolved by
the criterion of Legendre. In particular for a non-singular projective model $C$ we find
that $C(\mathbb{Q})$ is non-empty if and only if $C$ has p-adic points for all primes $p$, and this in
turn is determined by a finite number of congruences. If $C(\mathbb{Q})$ is non-empty then $C$ is
parametrized by rational functions and there are infinitely many rational points.

The most elusive case is that of genus 1. There may or may not be rational solutions and no
method is known for determining which is the case for any given curve. Moreover when
there are rational solutions there may or may not be infinitely many. If a non-singular
projective model $C$ has a rational point then $C(\mathbb{Q})$ has a natural structure as an abelian
group with this point as the identity element. In this case we call $C$ an elliptic curve over
$\mathbb{Q}$. In 1922 Mordell proved
that this group is finitely generated, thus fulfilling an implicit assumption of Poincar\'{e}.
If $C$ is an elliptic curve over $\mathbb{Q}$ then
$$
C(\mathbb{Q}) \simeq \mathbb{Z}^r \oplus C(\mathbb{Q})^{tors}
$$
for some integer $r \geq 0$, where $C(\mathbb{Q})^{tors}$ is a finite abelian group.
The integer $r$ is called the rank of $C$. It is zero if and only if $C(\mathbb{Q})$ is finite.

We can
find an affine model for an elliptic curve over Q in Weierstrass form
$$
C: y^2 = x^3 + ax + b
$$
with $a, b \in \mathbb{Z}$. We let $\Delta$ denote the discriminant of the cubic and set
$$
N_p := \sharp {(x,y)| y^2 �\equiv x^3 + ax + b  (p)}
$$
and $a_p := p - N_p$.
Then we can define the incomplete L-series of C (incomplete because we omit the Euler
factors for primes $p|2\Delta$) by
$$
L(C, s) := \Pi_{p\nshortmid 2\Delta} (1 - a_pp^{-s} + p^{1-2s})^{-1}.
$$
We view this as a function of the complex variables and this Euler product is then known
to converge for $Re(s) > 3/2$. A conjecture going back to Hasse
on 1952 predicted that $L(C, s)$ should have a holomorphic continuation as a
function of $s$ to the whole complex plane. This has now been proved.

Birch and Swinnerton-Dyer conjectured that the Taylor expansion of $L(C, s)$ at $s = 1$ has
the form
$L(C, s) = c(s - 1)^r +$ higher order terms
with $c\neq 0$ and $r = rank(C(\mathbb{Q}))$.
In particular this conjecture asserts that $L(C, 1) = 0$ �iff  $C(\mathbb{Q})$ is infinite.(Wiles 2006)

Early History: Problems on curves of genus 1 feature prominently in Diophantus' Arithmetica.
It is easy to see that a straight line meets an elliptic curve in three points (counting
multiplicity) so that if two of the points are rational then so is the third. In particular
if a tangent is taken to a rational point then it meets the curve again in a rational point.
Diophantus implicitly uses this method to obtain a second solution from a first. However
he does not iterate this process and it is Fermat who first realizes that one can sometimes
obtain infinitely many solutions in this way. Fermat also introduced a method of 'descent'
which sometimes permits one to show that the number of solutions is finite or even zero.
One very old problem concerned with rational points on elliptic curves is the congruent
number problem. One way of stating it is to ask which rational integers can occur as the
areas of right-angled triangles with rational length sides. Such integers are called congruent
numbers. For example, Fibonacci was challenged in the court of Frederic II with the
problem for $n = 5$ and he succeeded in finding such a triangle. He claimed moreover that
there was no such triangle for $n = 1$ but the proof was fallacious and the first correct proof
was given by Fermat. The problem dates back to Arab manuscripts of the 10th century. It is closely related to
the problem of determining the rational points on the curve $C_n: y^2 = x^3 - n^2x$. Indeed
$C_n(\mathbb{Q})$ is infinite iff $n$ is a congruent number.
Assuming the Birch and Swinnerton-Dyer conjecture (or even the weaker statement that
$C_n(\mathbb{Q})$ is infinite iff $L(C_n, 1) = 0$) one can show that any $n \equiv 5, 6, 7$ mod 8 is a congruent
number and moreover Tunnell has shown, again assuming the conjecture, that for n odd
and square-free
n is a congruent number iff
$$
\sharp\{x,y, z \in \mathbb{Z}: 2x^2 + y^2 + 8z^2 = n\}
= 2\times \sharp\{x, y, z \in \mathbb{Z}: 2x^2 + y^2 + 32z^2 = n\},
$$
with a similar criterion if n is even.

Recent History: It was the 1901 paper of Poincar\'{e} which started the modern interest in
the theory of rational points on curves and which first raised questions about the minimal
number of generators of $C(\mathbb{Q})$. The conjecture itself was first stated in the form we have
given in the early 1960's. In the intervening years the theory of $L$-functions of
elliptic curves (and other varieties) had been developed by a number of authors but the
conjecture was the first link between the $L$-function and the structure of $C(\mathbb{Q})$. It was
found experimentally using one of the early computers EDSAC at Cambridge. The first
general result proved was for elliptic curves with complex multiplication. (The curves with
complex multiplication fall into a finite number of families including ${y^2 = x^3 - Dx}$ and
${y^2 = x^3 - k}$ for varying $D, k
= 0$.) This theorem was proved in 1976 and is due to
Coates and Wiles. It states that if $C$ is a curve with complex multiplication and
$L(C, 1)= 0$ then $C(\mathbb{Q})$ is finite. In 1983 Gross and Zagier showed that if $C$ is a modular
elliptic curve and $L(C, 1) = 0$ but $L'(C, 1)\neq 0$, then an earlier construction of Heegner
actually gives a rational point of infinite order. Using new ideas together with this result,
Kolyvagin showed in 1990 that for modular elliptic curves, if $L(C, 1)\neq 0$ then $r = 0$ and if
$L(C, 1) = 0$ but $L'(C, 1)\neq 0$ then $r = 1$. In the former case Kolyvagin needed an analytic
hypothesis which was confirmed soon afterwards; Therefore if $L(C, s) \sim c(s - 1)^m$ with $c\neq 0$ and $m = 0$ or 1 then the conjecture holds.
In the cases where $m = 0$ or 1 some more precise results on c (which of course depends on
the curve) are known by work of Rubin and Kolyvagin.

\subsection{Hyperbolic algebraic varieties}

S. Lang in search for finiteness results, considers Kobayashi hyperbolicity, in which there is an
interplay between five notions:analytic notions of distance and measure;
complex analytic notions;
differential geometric notions of curvature (Chern and Ricci form);
algebraic notions of "general type" (pseudo ampleness);
arithmetic notions of rational points (existence of sections).
Lang is especially interested in the relations of the first four notions with
diophantine geometry, which historically has intermingled with complex differential
geometry. (S. Lang 1986)

Let $X$ denote a non-singular variety.
Among the possible complex analytic properties of $X$ we shall emphasize
that of being hyperbolic. There are several equivalent definitions of this notion,
and one of them, due to Brody, is that every holomorphic map of $C$ into $X$ is
constant.

The set of $F$-rational points of $X$ is denoted by $X(F)$. When we speak of
rational points, we shall always mean rational points in some field of finite type
over the rationals, that is, finitely generated over the rationals.

We shall say that $X$ (or a Zariski open subset) is mordellic if it has only a
finite number of rational points in every finitely generated field over $\mathbb{Q}$
according to our convention. Define the analytic exceptional set $Exc( X)$ of $X$ to
be the Zariski closure of the union of all images of non-constant holomorphic
maps $f:\mathbb{C} \to X$. Thus $X$ is hyperbolic if and only if this exceptional set is
empty. In general the exceptional set may be the whole variety, $X$ itself. Lang
conjectures that $X \setminus
 Exc( X)$ is mordellic.
It is also a problem to give an algebraic description of the exceptional set,
giving rise to the converse problem of showing that the exceptional set is not
mordellic, and in fact always has infinitely many rational points in a finite
extension of a field of definition for $X$. Define the algebraic exceptional set
$Exc_{alg}(X)$ to be the union of all non-constant rational images of $\mathbb{P}^1$ and
abelian varieties in $X$. Lang conjectures that the
analytic exceptional set is equal to the algebraic one. A subsidiary conjecture is
therefore that $Exc_{alg}(X)$ is closed, and that $X$ is hyperbolic if and only if every
rational map of $\mathbb{P}^1$ or an abelian variety into $X$ is constant. Similarly, until the
equality between the two exceptional sets is proved, one has the corresponding
conjecture that the complement of the algebraic exceptional set is mordellic.
Observe that the equality
$Exc(X) = Exc_{alg}(X)$
would give an algebraic characterization of hyperbolicity. Such a characterization
implies for instance that if a variety is defined over a field $F$ as above,
and is hyperbolic in one imbedding of $F$ in $\mathbb{C}$, then it is hyperbolic in every
imbedding of $F$ in $\mathbb{C}$, something which is by no means obvious, given the
analytic definitions of "hyperbolicity".
one wishes to characterize those varieties such that the exceptional set is a
proper subset. We shall give conjecturally a number of equivalent conditions,
which lead us into complex differential geometry, and algebraic geometry as
well as measure theory.

The properties having to do with hyperbolicity from the point of view of
differential geometry have been studied especially by Grauert-Reckziegel,
Green, Griffiths, and Kobayashi. Such properties have to do with "curvature".
We also consider a weakening of this notion,

namely measure hyperbolic, due to Kobayashi. The difference lies in looking at
positive $(1, l)$-forms or volume $(n, n)$-forms on a variety of dimension $n$. One
key aspect is that of the Ahlfors-Schwarz lemma, which states that a holomorphic
map is measure decreasing under certain conditions.
One associates a $(1, 1)$-form, the Ricci form $Ric(\Psi)$, to a
volume form $\Psi$ . The positivity of the Ricci form is related in a fundamental
way to hyperbolicity.

The problem is whether the following conditions are
equivalent:

1. Every sub variety of $X$ (including $X$ itself) is pseudo canonical.

2. $X$ is mordellic.

3. $X$ is hyperbolic.

3'. Every rational map of $\mathbb{P}^1$ or an abelian variety into $X$ is constant.

4. If $X$ is non-singular, there exist a positive $(1, l)$-form $\omega$ and a number
$B > 0$ such that for every complex one-dimensional immersed submanifold
$Y$ in $X$,
$Ric(\omega|Y)\geq B\omega|Y$.

5. Every subvariety of $X$ (including $X$ itself) is measure hyperbolic.

Evidence for the diophantine conjectures comes from their self-coherence,
rather than special cases.

Some of the listed conditions, like 1,2,3' are algebraic. The
others are analytic.

\section{Categorization of EPH phenomena}

The EPH-classifications of "Conic sections" and
"Action of $SL_2(\mathbb{R})$ on complex numbers, double numbers, and dual numbers" share appearance of
solid ellipses, parabolas and hyperbolas.
The EPH-classifications of "Elements of $SL_2(\mathbb{R})$" and
"Petersson's elliptic, parabolic, hyperbolic expansions of modular forms" match each other and can be related to the previous two
by considering $SO(2,1)$-representations which make solid ellipses, parabolas and hyperbolas appear in these cases also.  We call this the "conic sections" category. Limiting property of parabolic elements holds in all cases.

The EPH-classifications of "Euclidean and non-Euclidean geometries" and that of "Manifolds of constant sectional curvature"
and "Totally geodesic 2-dimensional foliations on 4-manifolds" and also "Riemann's uniformization theorem" share the concept of admitting a metric
of constant curvature. The EPH-classification of "Rational points on algebraic curves" also shares this feature but it is mainly
related to arithmetic not the concept of curvature. We call this the "curvature" category. One can say that limiting property of parabolic elements holds in the first four cases. The parabolic case is richer because these possess a group structure. The case of "Rational points on algebraic curves" shares this property. Although limiting property of parabolic case does not hold in the latter.

The EPH-classifications of "Hyperbolic, Parabolic and Elliptic dynamical systems" and that of "Partial differential equations"
share the appearance of different types of eigenvalues of linear maps
in certain order and therefore lie in a common category of EPH-classifications. We call this the "eigenvalues" category.

The EPH-classification of "Normal affine surfaces with $\mathbb{C}^*$-actions" is formulated in terms of fixed points of the
$\mathbb{C}^*$-actions. In some sense, this belongs to the category of "Hyperbolic, Parabolic and Elliptic dynamical systems" with complex time
and in other sense it belongs to the category of "Elements of $SL_2(\mathbb{R})$". It is difficult for me to decide which one is closer or if these two categories are the same.

Categorization of EPH-classifications as metric or conformal seems to be two independent groups as shown by nine Cayley-Klein geometries.
The EPH-classificastions "Euclidean and non-Euclidean geometries" and "Manifolds of constant sectional curvature"
together with "Totally geodesic 2-dimensional foliations on 4-manifolds" seem to be metric classifications.
The EPH-classificastions "Riemann's uniformization theorem" and "Partial differential equations" together with
 seem to be conformal classifications.
The EPH-classificastions "Elements of $SL_2(\mathbb{R})$",
"Action of $SL_2(\mathbb{R})$ on complex numbers, double numbers, and dual numbers"
"Petersson's elliptic, parabolic, hyperbolic expansions of modular forms",
"Rational points on algebraic curves" and
"Hyperbolic algebraic varieties" seem to belong to a different EPH-classification which we call Arithmetic classification.
The EPH-classificastions
"Hyperbolic, Parabolic and Elliptic dynamical systems" and "Normal affine surfaces with $\mathbb{C}^*$-actions"seem to be different from all above and we shall call it dynamical classification.

 It seems to us that metric classifications show both features of limiting property and richness of parabolic world.
But conformal, arithmetic and dynamical classifications only show richness of parabolic world. Therefore, we have at least two different categories
to study. If conformal and arithmetic and dynamical classifications come from the same paradigm, it means that we shall search for finiteness in
hyperbolic world even for conformal classifications. S. Lang's conjectures tend to promote the believe that conformal and arithmetic classifications are the same.

\section{Examples of EPH abuse}

At first glance, "Hyperbolic algebraic varieties" have no elliptic or parabolic counterparts. This may be because S. Lang was interested in finiteness results and deleted the elliptic or parabolic cases altogether. One can consider $Exc_{alg}^{ell}(X)$ the algebraic closure of the union of all images of $\mathbb{P}^1$ in $X$, and $Exc_{alg}^{par}(X)$ the algebraic closure of the union of all images of abelian varieties or $\mathbb{C}^*$ in $X$. The case $X=Exc_{alg}^{ell}(X)$ could be called elliptic. For example the projective space of any diemnsion is elliptic. The case $X=Exc_{alg}^{par}(X)$ could be called parabolic. In our point of view, defining hyperbolicity without considering the elliticity and parabolicity is not of theorization value, and is considered abuse of language. Since connections with the paradigm of EPH-classifications can not be made.

\section{EPH controversies}

A. Kisil explores the isometric action of
$G = SL_2(\mathbb{R})$ on the homogeneous spaces
$G/H$ where $H$ is one dimensional subgroup of $G$. There are only three such subgroups
up to conjugacy. The names elliptic, parabolic and hyperbolic come about from the shape of the equidistant
orbits which are ellipses, parabolas, hyperbolas respectively. Thus the Lobachevsky
geometry is elliptic (not hyperbolic) in this terminology. Therefore, there is a controversy here, if klein abused the name hyperbolic or not.
Or may be we shall come up with a more delicate conclusion and that is the symmetry of EPH-classifications with respect to ellipses and hyperbolas.
On one hand it seems a rather nice proposition, but does not match the association of hyperbolicity with finiteness. This association holds in
"Hyperbolic, Parabolic and Elliptic dynamical systems" and "Rational points on algebraic curves" and also in "Hyperbolic algebraic varieties".



\begin{thebibliography}{999999}

\bibitem[Ar-Vi]{Ar-Vi} V.Arauju, M. Viana, {Hyperbolic Dynamical Systems,} Math. Archive 0804.3192, 2008.

\bibitem[Ca-Gh]{Ca-Gh} G. Cairns, E. Ghys, {Totally geodesic foliations on 4-manifolds,} J. Diff. Geom. 23(1986) 241-254.

\bibitem[Cox]{Cox} H. S. M. Coxeter, {Non-Euclidean Geometry.} University of Toronto Press, 1998.

\bibitem[Cou-Hi] R. Courant, D. Hilbert, {Methods of Mathematical Physics, II,} New York: Wiley-Interscience, 1962.

\bibitem[Fl-Za]{Fl-Za} H. Flenner, M. Zaidenberg, {Normal affine surfaces with $\mathbb{C}^*$-action,} Osaka J. Math. 40(2003), 981-1009.

\bibitem[Ha-Ka]{Ha-Ka} B. Hasselblatt, A. Katok, {Principal Structures} in Handbook of Dynamical Systems Volume 1, Elsevier, 2002.

\bibitem[Im-O'S]{Im-O'S} O. Imamoghlu, C. O'Sullivan, {Parabolic, hyperbolic and elliptic Poincar\`{e} series,} math. archiv. 0806.4398 2008.

\bibitem[Ki]{Ki} A. Kisil, {Isometric action of $SL_2 (\mathbb{R})$ on homogenuous spaces,} math. archiv. 0810.0368, 2009.

\bibitem[Ki]{Ki} V. Kisil, {Starting with the group $SL_2 (\mathbb{R})$,} Notices of the AMS, 2007.

\bibitem[Kl]{Kl} F. Klein, {Uber die sogenannte nicht-Euklidische geometrie,} Gesammelte Math. Abh. I (1921), 254�305,
311�343, 344�350, 353�383.

\bibitem[La]{La} S. Lang, {Hyperbolic and Diophantine analysis,} Bull. AMS V. 14, No. 2, April 1986

\bibitem[Mc]{Mc} Alan S. McRae, {Clifford Algebras and Possible Kinematics,} in Symmetry, Integrability and Geometry: Methods and Applications SIGMA 3 (2007).

\bibitem[Pe]{Pe} Penrose R., {The road to reality,} Alfred A. Knopf, New York, 2005.

\bibitem[Ya]{Ya} Yaglom I.M., {A simple non-Euclidean geometry and its physical basis: an elementary account of Galilean geometry and the Galilean principle of relativity,} Heidelberg Science Library, translated from the Russian
by A. Shenitzer, with the editorial assistance of B. Gordon, Springer-Verlag, New York � Heidelberg, 1979.

\bibitem[Wi]{Wi} A. Wiles, {The Birch and Swinnerton-Dyer Conjecture,} 2006.



Sharif University of Technology, \\rastegar@sharif.ir

\end{thebibliography}
\end{document}